\documentclass[a4paper,11pt]{amsart}
\usepackage{graphicx}
\usepackage{mathptmx}
\usepackage{amsmath}
\usepackage{amssymb}
\usepackage{enumitem}
\usepackage{xcolor}

\newmuskip\pFqmuskip

\newcommand*\pFq[6][8]{%
  \begingroup 
  \pFqmuskip=#1mu\relax
  \mathcode`=\string"8000
  \begingroup\lccode`\~=`\,
  \lowercase{\endgroup\let~}\pFqcomma
  F^{#2}_{#3}{\left(\genfrac..{0pt}{}{#4}{#5}\bigg|#6\right)}%
  \endgroup
}
\newcommand{\pFqcomma}{\mskip\pFqmuskip}

\newtheorem{theorem}{Theorem}
\newtheorem{lemma}[theorem]{Lemma}

\newtheorem{proposition}[theorem]{Proposition}

\begin{document}

\title[A note on degenerate derangement polynomials and numbers]{A note on degenerate derangement polynomials and numbers}

\author{Taekyun  Kim}
\address{Department of Mathematics, Kwangwoon University, Seoul 139-701, Republic of Korea}
\email{tkkim@kw.ac.kr}

\author{Dae San  Kim}
\address{Department of Mathematics, Sogang University, Seoul 121-742, Republic of Korea}
\email{dskim@sogang.ac.kr}

\author{Hyunseok Lee}
\address{Department of Mathematics, Kwangwoon University, Seoul 139-701, Republic of Korea}
\email{luciasconstant@kw.ac.kr}

\author{Lee-Chae Jang}
\address{ Graduate School of Education, Kon-Kuk University,  Seoul 143-701, Republic of Korea}
\email{lcjang@konkuk.ac.kr}

\subjclass[2010]{11B83; 11B73; 60G50}
\keywords{degenerate derangement polynomials; degenerate gamma distribution; degenerate Fubini polynomials; fully degenerate Bell polynomials; degenerate Stirling numbers}
\maketitle

\begin{abstract}
In this paper, we study the degenerate derangement polynomials and numbers, investigate some properties of those polynomials and numbers and explore their connections with the degenerate gamma distributions. In more detail, we derive their explicit expressions, recurrence relations and some identities involving the degenerate derangement polynomials and numbers and other special polynomials and numbers, which include the fully degenerate Bell polynomials, the degenerate Fubini polynomials and the degenerate Stirling numbers of both kinds. We also show that those polynomials and numbers are connected with the moments of some variants of the degenerate gamma distributions.
\end{abstract}

\section{Introduction and preliminaries}
A derangement is a permutation with no fixed points. In other words, a derangement is a permutation of the elements of a set that leaves no elements in their original places. The number of derangements of a set of size $n$ is called the $n$-th derangement number and denoted by $d_n$. The first few terms of the derangement number sequence $\{d_{n}\}_{n=0}^{\infty}$ are $d_{0}=1,\ d_{1}=0,\ d_{2}=1,\ d_{3}=2,\ d_{4}=9,\dots$. It was Pierre R\'emonde de Motmort who initiated the study of counting derangements in 1708 (see [1]). \par
Carlitz was the first one who studied degenerate versions of some special polynomials and numers, namely the degenerate Bernoulli polynomials and numbers and degenerate Euler polynomials and numbers. In recent years, the study of various degenerate versions of some special polynomials and numbers regained the  interests of quite a few mathematicians and yielded many interesting arithmetical and combinatorial results. It is remarkable that the study of degenerate versions is not just limited to polynomials but can be extended to transcendental functions like gamma functions (see [9,14]).\par
The aim of this paper is to study the degenerate derangement polynomials, which are a degenerate version of the derangement polynomials. Here the derangement polynomials are a natural extension of the derangement numbers. In more detail, we derive their explicit expressions, recurrence relations and some identities involving those polynomials and numbers and other special polynomials and numbers, which include the fully degenerate Bell polynomials, the degenerate Fubini polynomials and the degerate Stirling numbers of both kinds. We also introduce the higher-order degenerate derangement polynomials. Then we explore the degenerate gamma distributions as a degenerate version of the gamma distributions and show that the moments of distributions coming from some variants of degenerate gamma distributions are related to the degenerate derangement polynomials or the degenerate derangement numbers or the higher-order degenerate derangement polynomials. \par
For the rest of this section, we recall the necessary facts about the degenerate derangement polynomials and numbers and the degenerate exponential functions. \par
As is well known, the generating function of the derangement numbers is given by
\begin{equation}
\frac{1}{1-t}e^{-t}=\sum_{n=0}^{\infty}d_{n}\frac{t^{n}}{n!},\quad(\mathrm{see}\ [1,3,4,8,12,13]).\label{1}
\end{equation}
From \eqref{1}, we note that
\begin{equation}
d_{n}=n!\sum_{i=0}^{n}\frac{(-1)^{i}}{i!},\quad(n\ge 0),\quad(\mathrm{see}\ [8,10,12,13]).\label{2}
\end{equation}
The derangement polynomials are defined by the generating function as
\begin{equation}
\frac{1}{1-t}e^{(x-1)t}=\sum_{n=0}^{\infty}d_{n}(x)\frac{t^{n}}{n!},\quad(\mathrm{see}\ [12,13]).\label{3}
\end{equation}
By \eqref{3}, we get
\begin{align}
d_{n}(x)\ &=\ \sum_{l=0}^{n}\binom{n}{l}d_{l}x^{n-l}\label{4}\\
&=\ n!\sum_{l=0}^{n}\frac{(x-1)^{l}}{l!},\quad(n\ge 0). \nonumber	
\end{align}
Clearly, we have $d_{n}(0)=d_{n}$. \par
For any nonzero real number $\lambda$, the degenerate exponential function is defined as
\begin{equation}
e_{\lambda}^{x}(t)=(1+\lambda t)^{\frac{x}{\lambda}}=\sum_{n=0}^{\infty}\frac{(x)_{n,\lambda}}{n!}t^{n},\quad(\mathrm{see}\ [2,5,9,11,16]),\label{5}
\end{equation}
where $(x)_{0,\lambda}=1,\ (x)_{n,\lambda}=x(x-\lambda)\cdots(x-(n-1)\lambda),\,\,(n \ge 1)$. \par
For brevity we denote $e_{\lambda}^{1}(t)$ by $e_{\lambda}(t)$. In this paper, we study the degenerate derangement polynomials which are derived from the degenerate exponential function. \par
From the definition of degenerate derangement polynomials, we investigate some properties and recurrence relations and new identities associated with special numbers and polynomials.
\section{Degenerate derangement polynomials}
In light of \eqref{3}, we may consider the degenerate derangement polynomials which are given by
\begin{equation}
\frac{1}{1-t}e^{x-1}_{\lambda}(t)=\sum_{n=0}^{\infty}d_{n,\lambda}(x)\frac{t^{n}}{n!}. \label{6}
\end{equation}
When $x=0$, $d_{n,\lambda}=d_{n,\lambda}(0)$ are called the degenerate derangement numbers. \par
From \eqref{5} and \eqref{6}, we get
\begin{align}
\sum_{n=0}^{\infty}d_{n,\lambda}(x)\frac{t^{n}}{n!}\ &=\ \sum_{l=0}^{\infty}t^{l}\sum_{m=0}^{\infty}(x-1)_{m,\lambda}\frac{t^{m}}{m!} \label{7} \\
&=\ \sum_{n=0}^{\infty}\bigg(n!\sum_{m=0}^{n}\frac{(x-1)_{m,\lambda}}{m!}\bigg)\frac{t^{n}}{n!}. \nonumber	
\end{align}
Comparing the coefficients on both sides of \eqref{7}, we obtain the following proposition.
\begin{proposition}
	For $n\ge 0$, we have
	\begin{displaymath}
		d_{n,\lambda}(x)=n!\sum_{l=0}^{n}\frac{(x-1)_{l,\lambda}}{l!}.
	\end{displaymath}
	In particular, for $x=0$, we obtain
	\begin{displaymath}
		d_{n,\lambda}=n!\sum_{l=0}^{n}\frac{(-1)_{l,\lambda}}{l!}.
	\end{displaymath}
\end{proposition}
Now, we observe that
\begin{equation}
e_{\lambda}^{x-1}(t)=1+\sum_{n=1}^{\infty}\bigg(d_{n,\lambda}(x)-nd_{n-1,\lambda}(x)\bigg)\frac{t^{n}}{n!}.\label{8}	
\end{equation}
From \eqref{5} and \eqref{8}, we have
\begin{equation}
(x-1)_{n,\lambda}=d_{n,\lambda}(x)-nd_{n-1,
\lambda}(x),\label{9}	
\end{equation}
and
\begin{displaymath}
	(-1)_{n,\lambda}=d_{n,\lambda}-nd_{n-1,\lambda},\quad(n \ge 1).
\end{displaymath}
In addition, by \eqref{6}, we get
\begin{equation}
d_{n,\lambda}(x)=\sum_{l=0}^{n}\binom{n}{l}d_{l,\lambda}(x)_{n-l,\lambda},\quad(n\ge 0). \label{10}
\end{equation}
Therefore, by \eqref{10}, we obtain the following theorem.
\begin{theorem}
The following identities hold true:
\begin{align*}
d_{n,\lambda}(x)=\sum_{l=0}^{n}\binom{n}{l}d_{l,\lambda}(x)_{n-l,\lambda},\quad(n \ge 0), \\
(x-1)_{n,\lambda}=d_{n,\lambda}(x)-nd_{n-1,\lambda}(x),\quad(n \ge 1), \\
(-1)_{n,\lambda}=d_{n,\lambda}-nd_{n-1,\lambda},\quad(n\ge 1),\quad(n \ge 1).
\end{align*}
\end{theorem}
Replacing $t$ by $1-e_{\lambda}(t)$ in \eqref{6}, we get
\begin{align}
e_{\lambda}^{x-1}\big(1-e_{\lambda}(t)\big)\ &=\ e_{\lambda}(t)\sum_{l=0}^{\infty}d_{l,\lambda}(x)\frac{1}{l!}\big(1-e_{\lambda}(t)\big)^{l}\label{11} \\
&=\ \sum_{m=0}^{\infty}\frac{(1)_{m,\lambda}}{m!}t^{m}\sum_{j=0}^{\infty}\sum_{l=0}^{j}(-1)^{l}d_{l,\lambda}(x)S_{2,\lambda}(j,l)\frac{t^{j}}{j!}\nonumber \\
&=\ \sum_{n=0}^{\infty}\bigg(\sum_{j=0}^{n}\sum_{l=0}^{j}\binom{n}{j}(1)_{n-j,\lambda}(-1)^{l}d_{l,\lambda}(x)S_{2,\lambda}(j,l)\bigg)\frac{t^{n}}{n!}.\nonumber
\end{align}
Here $S_{2,\lambda}(n,l),\ (n\ge l )$, are the degenerate Stirling numbers of the second kind given either by
\begin{displaymath}
(x)_{n,\lambda}=\sum_{l=0}^{n}S_{2,\lambda}(n,l)(x)_{l},\quad(n\ge 0),
\end{displaymath}
or by
\begin{displaymath}
\frac{1}{m!}(e_{\lambda}(t)-1)^m=\sum_{n=m}^{\infty} S_{2,\lambda}(n,m)\frac{t^n}{n!},\quad(m\ge 0), \quad(\mathrm{see}\ [7]),
\end{displaymath}
where $(x)_{0}=1$, $(x)_{n}=x(x-1)\cdots(x-n+1)$, $(n\ge 1)$.
Alternatively, \eqref{11} is also given by
\begin{align}
e_{\lambda}^{x-1}\big(1-e_{\lambda}(t)\big)\ &=\ \sum_{m=0}^{\infty}(x-1)_{m,\lambda}\frac{1}{m!}\big(1-e_{\lambda}(t)\big)^{m} \label{12} \\
&=\ \sum_{n=0}^{\infty}\bigg(\sum_{m=0}^{n}(x-1)_{m,\lambda}(-1)^{m}S_{2,\lambda}(n,m)\bigg)\frac{t^{n}}{n!}.\nonumber	
\end{align}
Therefore, by \eqref{11} and \eqref{12}, we obtain the following theorem.
\begin{theorem}
For $n\ge 0$, we have
\begin{displaymath}
\sum_{j=0}^{n}\sum_{l=0}^{j}\binom{n}{j}(1)_{n-j,\lambda}(-1)^{l}d_{l,\lambda}(x)S_{2,\lambda}(j,l)=\sum_{j=0}^{n}(x-1)_{j,\lambda}(-1)^{j}S_{2,\lambda}(n,j).
\end{displaymath}
\end{theorem}
Recently, the degenerate Fubini polynomials are introduced as
\begin{equation}
\frac{1}{1-y(e_{\lambda}(t)-1)}=\sum_{n=0}^{\infty}F_{n,\lambda}(y)\frac{t^{n}}{n!},\quad(\mathrm{see}\ [11,15]).\label{13}	
\end{equation}
Note that $\displaystyle\lim_{\lambda\rightarrow 0}F_{n,\lambda}(y)=F_{n}(y)\displaystyle$ are the ordinary Fubini polynomials (see [6]). Replacing $t$ by $e_{\lambda}(t)-1$ in \eqref{6}, we get
\begin{align}
\frac{1}{2-e_{\lambda}(t)}e_{\lambda}^{x-1}\big(e_{\lambda}(t)-1\big)\ &=\ \sum_{l=0}^{\infty}d_{l,\lambda}(x)\frac{1}{l!}\big(e_{\lambda}(t)-1\big)^{l}\label{14}\\
&=\ \sum_{n=0}^{\infty}\bigg(\sum_{l=0}^{n}S_{2,\lambda}(n,l)d_{l,\lambda}(x)\bigg)\frac{t^{n}}{n!}.\nonumber
\end{align}
In terms of \eqref{13}, we note that \eqref{14} is also given by
\begin{align}
&\frac{1}{2-e_{\lambda}(t)}e_{\lambda}^{x-1}\big(e_{\lambda}(t)-1\big)\label{15}\\
&\ =\ \sum_{l=0}^{\infty}F_{l,\lambda}(1)\frac{t^{l}}{l!}\sum_{m=0}^{\infty}(x-1)_{m,\lambda}\frac{1}{m!}\big(e_{\lambda}(t)-1\big)^{m}\nonumber\\
&\ =\  \sum_{l=0}^{\infty}F_{l,\lambda}(1)\frac{t^{l}}{l!}\sum_{j=0}^{\infty}\sum_{m=0}^{j}(x-1)_{m,\lambda}S_{2,\lambda}(j,m)\frac{t^{j}}{j!}\nonumber \\
&\ =\ \sum_{n=0}^{\infty}\bigg(\sum_{l=0}^{n}\sum_{m=0}^{l}\binom{n}{l}F_{n-l,\lambda}(1)(x-1)_{m,\lambda}S_{2,\lambda}(l,m)\bigg)\frac{t^{n}}{n!}.\nonumber
\end{align}
Therefore, by \eqref{14} and \eqref{15}, we obtain the following theorem.
\begin{theorem}
	For $n\ge 0$, we have
	\begin{displaymath}
		\sum_{l=0}^{n}S_{2,\lambda}(n,l)d_{l,\lambda}(x)=\sum_{l=0}^{n}\sum_{m=0}^{l}\binom{n}{l}F_{n-l,\lambda}(1)(x-1)_{m,\lambda}S_{2,\lambda}(l,m).
	\end{displaymath}
\end{theorem}
Let $\log_{\lambda}(t)$ be the compositional inverse function of $e_{\lambda}(t)$. Recall that the degenerate Stirling numbers of the first kind are defined either by
\begin{displaymath}
(x)_{n}=\sum_{l=0}^{n}S_{1,\lambda}(n,l)(x)_{l,\lambda},\quad(n\ge 0),
\end{displaymath}
or by
\begin{displaymath}
\frac{1}{m!}\big(\log_{\lambda}(1+t)\big)^m=\sum_{n=m}^{\infty}S_{1,\lambda}(n,m)\frac{t^n}{n!},\quad(m \ge 0), \quad(\mathrm{see}\ [7,14]).
\end{displaymath} \par
Replacing $t$ by $\log_{\lambda}(1+t)$ in \eqref{13} with $y=1$, we get
\begin{align}
\frac{1}{1-t}\ &=\bigg(\frac{1}{1-t}e^{-t}\bigg)e^{t}\label{17} \\
&=\ \sum_{l=0}^{\infty}F_{l,\lambda}(1)\frac{1}{l!}\big(\log_{\lambda}(1+t)\big)^{l}\nonumber \\	
&=\ \sum_{n=0}^{\infty}\bigg(\sum_{l=0}^{n}F_{l,\lambda}(1)S_{1,\lambda}(n,l)\bigg)\frac{t^{n}}{n!}.\nonumber
\end{align}
Writing the left hand side of \eqref{17} differently, we have
\begin{align}
\frac{1}{1-t}\ &=\ \bigg(\frac{1}{1-t}e^{-1}_{\lambda}(t)\bigg)e_{\lambda}(t) \label{18} \\ &=\sum_{l=0}^{\infty}d_{l,\lambda}\frac{t^{l}}{l!}\sum_{m=0}^{\infty}(1)_{m,\lambda}\frac{t^{m}}{m!}\nonumber \\
&=\ \sum_{n=0}^{\infty}\bigg(\sum_{l=0}^{n}\binom{n}{l}d_{l,\lambda}(1)_{n-l,\lambda}\bigg)\frac{t^{n}}{n!}.\nonumber	
\end{align}
Therefore, by \eqref{6}, \eqref{17} and \eqref{18}, we obtain the following theorem.
\begin{theorem}
	For $n\ge 0$, we have
	\begin{displaymath}
		\sum_{l=0}^{n}F_{l,\lambda}(1)S_{1,\lambda}(n,l)=\sum_{l=0}^{n}\binom{n}{l}d_{l,\lambda}(1)_{n-l,\lambda}=\sum_{l=0}^{n}\binom{n}{l}d_{l,\lambda}(x)(1-x)_{n-l,\lambda}.
	\end{displaymath}
\end{theorem}
Replacing $t$ by $e_{\lambda}(t)-1$ in \eqref{8}, we obtain
\begin{align}
e_{\lambda}^{x-1}\big(e_{\lambda}(t)-1\big)\ &=\ 1+\sum_{m=1}^{\infty}\big(d_{m,\lambda}(x)-md_{m-1,\lambda}(x)\big)\frac{1}{m!}\big(e_{\lambda}(t)-1\big)^{m}\label{19}\\
&=\ 1+\sum_{m=1}^{\infty}\big(d_{m,\lambda}(x)-md_{m-1,\lambda}(x)\big)\sum_{n=m}^{\infty}S_{2,\lambda}(n,m)\frac{t^{n}}{n!} \nonumber \\
&=\ 1+\sum_{n=1}^{\infty}\bigg(\sum_{m=1}^{\infty}\big(d_{m,\lambda}(x)-md_{m-1,\lambda}(x)\big)S_{2,\lambda}(n,m)\bigg)\frac{t^{n}}{n!}. \nonumber
\end{align}
We note that \eqref{19} is also given by
\begin{align}
e_{\lambda}^{x-1}\big(e_{\lambda}(t)-1\big)\ &=\ \sum_{m=0}^{\infty}(x-1)_{m,\lambda}\frac{1}{m!}\big(e_{\lambda}(t)-1
\big)^{m}\label{20}\\
&=\ \sum_{m=0}^{\infty}(x-1)_{m,\lambda}\sum_{n=m}^{\infty}S_{2,\lambda}(n,m)\frac{t^{n}}{n!} \nonumber \\
&=\ 1+\sum_{n=1}^{\infty}\bigg(\sum_{m=1}^{n}(x-1)_{m,\lambda}S_{2,\lambda}(n,m)\bigg)\frac{t^{n}}{n!}. \nonumber
\end{align}
Therefore, by \eqref{19} and \eqref{20}, we obtain the following lemma.
\begin{lemma}
	For $n\in\mathbb{N}$, we have
	\begin{displaymath}
		\sum_{m=1}^{n}(x-1)_{m,\lambda}S_{2,\lambda}(n,m)=\sum_{m=1}^{n}\big(d_{m,\lambda}(x)-md_{m-1,\lambda}(x)\big)S_{2,\lambda}(n,m).
	\end{displaymath}
	In particular, for $x=0$, we get
	\begin{displaymath}
		\sum_{m=1}^{n}(-1)_{m,\lambda}S_{2,\lambda}(n,m)=\sum_{m=1}^{n}\big(d_{m,\lambda}-md_{m-1,\lambda}\big)S_{2,\lambda}(n,m).
	\end{displaymath}
\end{lemma}
Recently, Kim-Kim considered the fully degenerate Bell polynomials given by
\begin{equation}
e_{\lambda}\big(x(e_{\lambda}(t)-1)\big)=\sum_{n=0}^{\infty}\mathrm{Bel}_{n,\lambda}(x)\frac{t^{n}}{n!},\quad(\mathrm{see}\ [2]).\label{21}
\end{equation}
Replacing $t$ by $\log_{\lambda}(1+t)$ in \eqref{21} with $x=1$, we get
\begin{align}
e_{\lambda}(t)\ &=\ \sum_{m=0}^{\infty}\mathrm{Bel}_{m,\lambda}\frac{1}{m!}\big(\log_{\lambda}(1+t)\big)^{m}\label{22} \\
&=\ \sum_{m=0}^{\infty}\mathrm{Bel}_{m,\lambda}\sum_{n=m}^{\infty}S_{1,\lambda}(n,m)\frac{t^{n}}{n!} \nonumber\\
&=\ \sum_{n=0}^{\infty}\bigg(\sum_{m=0}^{n}\mathrm{Bel}_{m,\lambda}S_{1,\lambda}(n,m)\bigg)\frac{t^{n}}{n!}.\nonumber
\end{align}
Obviously, \eqref{22} is also given by
\begin{equation}
e_{\lambda}(t)=\sum_{n=0}^{\infty}(1)_{n,\lambda}\frac{t^{n}}{n!}. \label{23}	
\end{equation}
Therefore, by \eqref{22} and \eqref{23}, we obtain the following theorem.
\begin{theorem}
	For $n\ge 0$, we have
	\begin{displaymath}
		(1)_{n,\lambda}=\sum_{m=0}^{n}\mathrm{Bel}_{m,\lambda}S_{1,\lambda}(n,m),
	\end{displaymath}
	and
	\begin{displaymath}
		\mathrm{Bel}_{n,\lambda}=\sum_{m=0}^{n}(1)_{m,\lambda}S_{2,\lambda}(n,m).
	\end{displaymath}
\end{theorem}
Now, we observe from \eqref{6} that
\begin{align}
\frac{1}{1-t}\ &=\bigg(\ \sum_{m=0}^{\infty}d_{m,\lambda}(x)\frac{t^{m}}{m!}\bigg) e_{\lambda}^{1-x}(t) \label{24}\\
&=\ \sum_{m=0}^{\infty}d_{m,\lambda}(x)\frac{t^{m}}{m!}\sum_{l=0}^{\infty}(1-x)_{l,\lambda}\frac{t^{l}}{l!} \nonumber\\
&=\ \sum_{n=0}^{\infty}\bigg(\sum_{m=0}^{n}\binom{n}{m}d_{m,\lambda}(x)(1-x)_{n-m,\lambda}\bigg)\frac{t^{n}}{n!}.\nonumber	
\end{align}
Another expression of \eqref{24} is given by
\begin{align}
\frac{1}{1-t}\ &=\ e_{\lambda}^{-1}\big(\log_{\lambda}(1-t)\big)\ =\ \sum_{m=0}^{\infty}(-1)_{m,\lambda}\frac{1}{m!}\big(\log_{\lambda}(1-t)\big)^{m}\label{25} \\
&=\ \sum_{m=0}^{\infty}	(-1)_{m,\lambda}\sum_{n=m}^{\infty}(-1)^{n}S_{1,\lambda}(n,m)\frac{t^{n}}{n!}=\sum_{n=0}^{\infty}\bigg(\sum_{m=0}^{n}(-1)_{m,\lambda}(-1)^{n}S_{1,\lambda}(n,m)\bigg)\frac{t^{n}}{n!}. \nonumber
\end{align}
From \eqref{24} and \eqref{25}, we note that
\begin{displaymath}
	(-1)^{n}\sum_{m=0}^{n}(-1)_{m,\lambda}S_{1,\lambda}(n,m)=\sum_{m=0}^{n}\binom{n}{m}d_{m,\lambda}(x)(1-x)_{n-m,\lambda}.
\end{displaymath} \par
Replacing $t$ by $\log_{\lambda}(1-t)$ in \eqref{21} with $x=1$, we get
\begin{align}
e_{\lambda}(-t)\ &=\ \sum_{k=0}^{\infty}\mathrm{Bel}_{k,\lambda}\frac{1}{k!}\big(\log_{\lambda}(1-t)\big)^{k}\nonumber\\
&=\ \sum_{k=0}^{\infty}\mathrm{Bel}_{k,\lambda}\sum_{n=k}^{\infty}S_{1,\lambda}(n,k)(-1)^{n}\frac{t^{n}}{n!}\label{26}\\
&=\ \sum_{n=0}^{\infty}\bigg(\sum_{k=0}^{n}\mathrm{Bel}_{k,\lambda}S_{1,\lambda}(n,k)(-1)^{n}\bigg)\frac{t^{n}}{n!}.\nonumber	
\end{align}
We remark that \eqref{26} is alternatively given by
\begin{equation}
e_{\lambda}(-t) \ =\ e_{-\lambda}^{-1}(t)\ =\ \sum
_{n=0}^{\infty}(-1)_{n,-\lambda}\frac{t^n}{n!}. \label{27} \\
\end{equation}
Thus, from \eqref{26} and \eqref{27}, we have
\begin{equation}
\sum_{k=0}^{n}\mathrm{Bel}_{k,\lambda}S_{1,\lambda}(n,k)\ =\ (-1)^{n}(-1)_{n,-\lambda},\quad(n\ge 0).\label{28}
\end{equation} \par
Replacing $t$ by $1-e_{-\lambda}(t)$ in \eqref{6} with $x=0$, we get
\begin{align}
e_{-\lambda}^{-1}(t)e_{\lambda}^{-1}\big(1-e_{-\lambda}(t)\big)\ &=\ \sum_{m=0}^{\infty}d_{m,\lambda}\frac{(-1)^{m}}{m!}\big(e_{-\lambda}(t)-1\big)^{m} \label{29}\\
&=\ \sum_{m=0}^{\infty}d_{m,\lambda}(-1)^{m}\sum_{n=m}^{\infty}S_{2,-\lambda}(n,m)\frac{t^{n}}{n!} \nonumber \\
&=\ \sum_{n=0}^{\infty}\bigg(\sum_{m=0}^{n}d_{m,\lambda}(-1)^{m}S_{2,-\lambda}(n,m)\bigg)\frac{t^{n}}{n!}.\nonumber
\end{align}
An alternative expression of \eqref{29} is given by
\begin{align}
e_{-\lambda}^{-1}(t)e_{\lambda}^{-1}\big(1-e_{-\lambda}(t)\big)\ &=\ e_{-\lambda}^{-1}(t) e_{-\lambda}\big(e_{-\lambda}(t)-1\big) \label{30} \\
&=\ \sum_{l=0}^{\infty}(-1)_{l,-\lambda}\frac{t^{l}}{l!}\sum_{m=0}^{\infty}\mathrm{Bel}_{m,-\lambda}\frac{t^{m}}{m!} \nonumber\\
&=\ \sum_{n=0}^{\infty}\bigg(\sum_{m=0}^{n}\binom{n}{m}\mathrm{Bel}_{m,-\lambda}(-1)_{n-m,-\lambda}\bigg)\frac{t^{n}}{n!}. \nonumber	
\end{align}
From \eqref{29} and \eqref{30}, we have
\begin{equation}
\sum_{m=0}^{n}(-1)^{m}d_{m,\lambda}S_{2,-\lambda}(n,m)=\sum_{m=0}^{n}\binom{n}{m}\mathrm{Bel}_{m,-\lambda}(-1)_{n-m,-\lambda},\quad(n\ge 0).\label{31}
\end{equation}
Therefore, by \eqref{28} and \eqref{31}, we obtain the following theorem.
\begin{theorem}
For $n\ge 0$, we have
\begin{displaymath}
\sum_{m=0}^{n}(-1)^{m}d_{m,\lambda}S_{2,-\lambda}(n,m)=\sum_{m=0}^{n}\binom{n}{m}\mathrm{Bel}_{m,-\lambda}(-1)_{n-m,-\lambda}.
\end{displaymath}
In addition, we have
\begin{displaymath}
\sum_{k=0}^{n}\mathrm{Bel}_{k,\lambda}S_{1,\lambda}(n,k)=(-1)^{n}(-1)_{n,-\lambda},\quad(n \ge 0).
\end{displaymath}
\end{theorem}
For $r\in\mathbb{N}$, we define the degenerate derangement polynomials of order $r$ which are given by
\begin{equation}
\frac{1}{(1-t)^{r}}e_{\lambda}^{x-1}(t)=\sum_{n=0}^{\infty}d_{n}^{(r)}(x)\frac{t^{n}}{n!}. \label{32}
\end{equation}
When $x=0$, $d_{n}^{(r)}(0)$ are called the degenerate derangement numbers of order $r$. \par
~~\\
From \eqref{32}, we note that
\begin{align}
\sum_{n=0}^{\infty}d_{n}^{(r)}(x)\frac{t^{n}}{n!}\ &=\ \sum_{m=0}^{\infty}\binom{r+m-1}{m}t^{m}\sum_{l=0}^{\infty}(x-1)_{l,\lambda}\frac{t^{l}}{l!}\label{33}\\
&=\ \sum_{n=0}^{\infty}\bigg(n!\sum_{l=0}^{n}\frac{(x-1)_{l,\lambda}}{l!}\binom{r+n-l-1}{n-l}\bigg)\frac{t^{n}}{n!}.\nonumber
\end{align}
Comparing the coefficients on both sides of \eqref{33}, we obtain the following theorem.
\begin{theorem}
	For $n\ge 0$, we have
	\begin{displaymath}
		d_{n}^{(r)}(x)=n!\sum_{l=0}^{n}\frac{(x-1)_{l,\lambda}}{l!}\binom{r+n-l-1}{n-l}.
	\end{displaymath}
	In particular, for $x=0$, we have
	\begin{displaymath}
		d_{n}^{(r)}=n!\sum_{l=0}^{n}\frac{(-1)_{l,\lambda}}{l!}\binom{r+n-l-1}{n-l}.
	\end{displaymath}
\end{theorem}
By \eqref{6}, we get
\begin{equation}
\frac{1}{1+t}e_{\lambda}^{-1}(-t)=\sum_{m=0}^{\infty}d_{m,\lambda}(-1)^{m}\frac{t^{m}}{m!}. \label{34}
\end{equation}
Replacing $t$ by $e_{-\lambda}(t)-1$ in $\eqref{34}$, we get
\begin{align}
e_{\lambda}^{-1}\big(1-e_{-\lambda}(t)\big)\ &=\ e_{-\lambda}(t)\sum_{m=0}^{\infty}d_{m,\lambda}(-1)^{m}\frac{1}{m!}\big(e_{-\lambda}(t)-1\big)^{m}\label{35} \\
&=\ e_{-\lambda}(t)\sum_{m=0}^{\infty}d_{m,\lambda}(-1)^{m}\sum_{j=m}^{\infty}S_{2,-\lambda}(j,m)\frac{t^{j}}{j!}\nonumber \\
&=\ \sum_{l=0}^{\infty}(1)_{l,-\lambda}\frac{t^{l}}{l!}\sum_{j=0}^{\infty}\bigg(\sum_{m=0}^{j}(-1)^{m}d_{m,\lambda}S_{2,-\lambda}(j,m)\bigg)\frac{t^{j}}{j!} \nonumber \\
&=\ \sum_{n=0}^{\infty}\bigg(\sum_{j=0}^{n}\sum_{m=0}^{j}\binom{n}{j}(1)_{n-j,-\lambda}(-1)^{m}d_{m,\lambda}S_{2,-\lambda}(j,m)\bigg)\frac{t^{n}}{n!}.\nonumber
\end{align}
Alternatively, \eqref{35} is also given by
\begin{equation}
e_{\lambda}^{-1}\big(1-e_{-\lambda}(t)\big)=e_{-\lambda}\big(e_{-\lambda}(t)-1\big)=\sum_{n=0}^{\infty}\mathrm{Bel}_{n,-\lambda}(1)\frac{t^{n}}{n!}. \label{36}
\end{equation}
Therefore, by \eqref{35} and \eqref{36}, we obtain the following theorem.
\begin{theorem}
	For $n\ge 0$, we have
	\begin{displaymath}
		\mathrm{Bel}_{n,-\lambda}(1)=\sum_{j=0}^{n}\sum_{m=0}^{j}\binom{n}{j}(1)_{n-j,-\lambda}(-1)^{m}d_{m,\lambda}S_{2,-\lambda}(j,m).
	\end{displaymath}
\end{theorem}
\section{Further remarks}
Let $f(x)$ be the probability density function of the continuous random variable $X$, and let $g(x)$ be a real valued function. Then the expectation of $g(X)$, $E[g(X)]$, is defined by
\begin{equation}
E[g(X)]=\int_{-\infty}^{\infty}g(x)f(x)dx,\quad(\mathrm{see}\ [18]).\label{37}
\end{equation}
A continuous random variable $X$, whose density function is given by
\begin{equation}
f(x)=\left\{\begin{array}{ccc}
	\beta e^{-\beta x}\frac{(\beta x)^{\alpha-1}}{\Gamma(\alpha)}, & \textrm{if $x\ge 0$}, \\
	0, & \textrm{ if $x<0$,}
\end{array}\right.\label{38}	
\end{equation}
for some $\beta>0$ and $\alpha>0$, is said to be the gamma random variable with parameters $\alpha,\beta$ and denoted by $X\sim\Gamma(\alpha,\beta)$. \par
Let $X\sim\Gamma(1,1)$. Then, for all $t < 1$, we have
\begin{align}
E\big[e^{Xt}\cdot e_{\lambda}^{-1}(t)\big]\ &=\ e_{\lambda}^{-1}(t)\int_{0}^{\infty}e^{xt}e^{-x}dx \nonumber\\
&=\ \frac{1}{1-t}e_{\lambda}^{-1}(t)\ =\ \sum_{n=0}^{\infty}d_{n,\lambda}\frac{t^{n}}{n!}.\label{39}
\end{align}
Clearly, we also have
\begin{equation}
E\big[e^{Xt} e_{\lambda}^{-1}(t)\big]=\sum_{n=0}^{\infty}\bigg(\sum_{m=0}^{n}\binom{n}{m}(-1)_{n-m,\lambda}E[X^{m}]\bigg)\frac{t^{n}}{n!}. \label{40}
\end{equation}
Therefore, by \eqref{39} and \eqref{40}, we obtain the following equations. \par
For $n\ge 0$, we have
\begin{displaymath}
\sum_{m=0}^{n}\binom{n}{m}(-1)_{n-m,\lambda}E[X^{m}]=d_{n,\lambda},
\end{displaymath}
and, more generally, we also have
\begin{displaymath}
\sum_{m=0}^{n}\binom{n}{m}(x-1)_{n-m,\lambda}E[X^{m}]=d_{n,\lambda}(x).
\end{displaymath} \par
Unless otherwise stated, for the rest of this section we assume that $\lambda \in (0,1)$.
The degenerate gamma function $\Gamma_{\lambda}(x)$, which is initially defined for $0<\mathrm{Re}(s)<\frac{1}{\lambda}$ by the following integral
\begin{equation}
\Gamma_{\lambda}(s)=\int_{0}^{\infty}e_{\lambda}^{-1}(t)t^{s-1}dt,\quad(\mathrm{see}\ [9,14]),\label{41}
\end{equation}
can be continued to a meromorphic function on $\mathbb{C}$, whose only singularities are simple poles at $s=0,-1,-2,\dots,\frac{1}{\lambda},\frac{1}{\lambda}+1,\frac{1}{\lambda}+2,\dots$.
Thus, by \eqref{41}, we get
\begin{equation}
\Gamma_{\lambda}(k)=\frac{\Gamma(k)}{(1)_{k+1,\lambda}},\quad\bigg(k\in\mathbb{N},\ \lambda \in (0,\frac{1}{k})\bigg),\label{42}
\end{equation}
and, in particular, we have
\begin{displaymath}
	\Gamma_{\lambda}(1)=\frac{1}{1-\lambda},\quad(\mathrm{see}\ [9]).
\end{displaymath}
A random variable $X=X_{\lambda}$ is said to have the degenerate gamma distribution with parameters $\alpha$ and $\beta$, $\big(\frac{1}{\lambda}>\alpha>0,\ \beta>0\big)$, and denoted by $X\sim\Gamma_{\lambda}(\alpha,\beta)$, if its probability density function has the form
\begin{displaymath}
f_{\lambda}(x)=\left\{\begin{array}{ccc}
\frac{1}{\Gamma_{\lambda}(\alpha)}\beta(\beta x)^{\alpha-1}e_{\lambda}^{-1}(\beta x), & \textrm{if $x \ge 0$,}\\
0, & \textrm{otherwise.}
\end{array}\right.
\end{displaymath} \par
Note that $\frac{d}{dx}e_{\lambda}^{c}(x)=ce_{\lambda}^{c-\lambda}(x)$, for any constant $c$. Then, for $X\sim\Gamma_{\lambda}(1,1)$, we have
\begin{align}
E\big[e_{\lambda}^{t-\lambda}(X)\big]\ &=\ (1-\lambda)\int_{0}^{\infty}e_{\lambda}^{t-\lambda}(x)e_{\lambda}^{-1}(x)dx	 \label{43}\\
&=\ (1-\lambda)\int_{0}^{\infty}e_{\lambda}^{t-1-\lambda}(x)dx\ =\ \frac{1}{1-\lambda}\frac{1}{1-t}e_{\lambda}^{-1}(t)e_{\lambda}(t)\nonumber\\
&=\ (1-\lambda)\sum_{l=0}^{\infty}d_{l,\lambda}\frac{t^{l}}{l!}\sum_{m=0}^{\infty}(1)_{m,\lambda}\frac{t^{m}}{m!}\nonumber\\
&=\ \sum_{n=0}^{\infty}(1-\lambda)\sum_{l=0}^{n}d_{l,\lambda}(1)_{n-l,\lambda}\binom{n}{l}\frac{t^{n}}{n!}.\nonumber
\end{align}
Evidently, we also have
\begin{align}
E\big[e_{\lambda}^{t-\lambda}(X)\big]\ &=\ E\bigg[\frac{1}{1+\lambda X}\big(1+\lambda X\big)^{\frac{t}{\lambda}}\bigg]\nonumber \\
&=\ \sum_{n=0}^{\infty}E\bigg[\frac{1}{1+\lambda X}\bigg(\frac{1}{\lambda}\log(1+\lambda X)\bigg)^{n}\bigg]\frac{t^{n}}{n!}. \label{44}
\end{align}
Therefore, \eqref{43} and \eqref{44}, we obtain the following theorem.
\begin{theorem}
For $X\sim\Gamma_{\lambda}(1,1)$, we have
\begin{displaymath}
E\bigg[\frac{1}{1+\lambda X}\bigg(\frac{1}{\lambda}\log(1+\lambda X)\bigg)^{n}\bigg]=(1-\lambda)\sum_{l=0}^{n}d_{l,\lambda}(1)_{n-l,\lambda}\binom{n}{l}.
\end{displaymath}
\end{theorem}
Now, we observe that
\begin{displaymath}
	\big(\log(1+\lambda X)\big)^{n}=n!\sum_{m=n}^{\infty}S_{1}(m,n)\frac{\lambda^{m}}{m!}X^{m},\quad (n\ge 0),
\end{displaymath}
where $S_{1}(n,m)$ are the Stirling numbers of the first kind, (see [17,19,20]).
In turn, we have
\begin{equation}
E\bigg[\frac{1}{1+\lambda X}\bigg(\frac{1}{\lambda}\log(1+\lambda X)\bigg)^{n}\bigg]=
\frac{n!}{\lambda^n}\sum_{m=n}^{\infty}S_{1}(m,n)\frac{\lambda^{m}}{m!}E\bigg[\frac{X^{m}}{1+\lambda X}\bigg].\label{45}
\end{equation}
From Theorem 11 and \eqref{45}, we have
\begin{displaymath}
\sum_{n=m}^{\infty}S_{1}(n,m)\frac{\lambda^{m}}{m!}E\bigg[\frac{X^{m}}{1+\lambda X}\bigg]=(1-\lambda)\frac{\lambda^n}{n!}\sum_{l=0}^{n}d_{l,\lambda}(1)_{n-l,\lambda}\binom{n}{l},\quad(n\ge 0),
\end{displaymath}
where $X\sim\Gamma_{\lambda}(1,1)$. \par
For $X_{1},X_{2},\dots,X_{r}\sim\Gamma(1,1)$, assume that $X_{1},X_{2},\dots,X_{r}$ are independent. Then we have
\begin{align}
E\big[e^{(X_{1}+X_{2}+\cdots+X_{r})t}e_{\lambda}^{x-1}(t)\big]\ &=\ E\big[e^{X_{1}t}\big] E\big[e^{X_{2}t}\big]\cdots E\big[e^{X_{r}t}\big]\cdot e_{\lambda}^{x-1}(t) \label{46} \\
&=\ \underbrace{\bigg(\frac{1}{1-t}\bigg)\times\bigg(\frac{1}{1-t}\bigg)\times\cdots\times\bigg(\frac{1}{1-t}\bigg)}_{r-\mathrm{times}}e_{\lambda}^{x-1}(t) \nonumber\\
&=\ \sum_{n=0}^{\infty}d_{n}^{(r)}(x)\frac{t^{n}}{n!}.\nonumber
\end{align}
Alternatively, \eqref{46} is given by
\begin{align}
&E\big[e^{(X_{1}+\cdots+X_{r})t}e_{\lambda}^{x-1}(t)\big]\label{47}\\
 &\quad =\ \sum_{l=0}^{\infty}E\big[(X_{1}+\cdots+X_{r})^{l}\big]\frac{t^{l}}{l!}\sum_{m=0}^{\infty}(x-1)_{m,\lambda}\frac{t^{m}}{m!}\nonumber \\
 &\quad =\ \sum_{n=0}^{\infty}\bigg(\sum_{l=0}^{n}\binom{n}{l}E\big[(X_{1}+\cdots+X_{r})^{l}\big](x-1)_{n-l,\lambda}\bigg)\frac{t^{n}}{n!}. \nonumber
\end{align}
By \eqref{46} and \eqref{47}, we get
\begin{displaymath}
	d_{n,\lambda}^{(r)}(x)=\sum_{l=0}^{n}\binom{n}{l}E\big[(X_{1}+\cdots+X_{r})^{l}\big](x-1)_{n-l,\lambda},\quad(n\ge 0).
\end{displaymath}

\section{Conclusion}

In this paper, we have dealt with the degenerate derangement polynomials $d_{n,\lambda}(x)$, which are a degenerate version of the derangement polynomials $d_n(x)$. We derived their explicit expressions, recurrence relations and some identities involving those polynomials and numbers and other special polynomials and numbers such as the fully degenerate Bell polynomials, the degenerate Fubini polynomials and the degerate Stirling numbers of both kinds. We also introduced the higher-order degenerate derangement polynomials. Then we explored the degenerate gamma distributions as a degenerate version of the gamma distributions and showed that the moments of distributions coming from some variants of degenerate gamma distributions are related to the degenerate derangement polynomials or the degenerate derangement numbers or the higher-order degenerate derangement polynomials.\par
In recent years, the study of many special numbers and polynomials has been carried out by using several different methods, which include generating functions, combinatorial methods, umbral calculus, $p$-adic analysis, probability theory, special functions and differential equations. Moreover, the same has been done for various degenerate versions of quite a few special numbers and polynomials. Motivations for studying degenerate versions arise from their interests not only in combinatorial and arithmetical properties but also in their applications to symmetric identities, differential equations and probability theories. \par
It is one of our future projects to continue to investigate many ordinary and degenerate special numbers and polynomials by various means and to find their applications in physics, science and engineering as well as in mathematics.

\end{document}